\input amstex
\magnification 1200
\documentstyle{amsppt}
\vcorrection{-1cm}

\def\refAS    {1}
\def\refBZ    {2}
\def\refLin   {3}
\def\refLNM   {4}
\def\refLNMii {5}
\def\refM     {6}
\def\refPalka {7}
\def\refTono  {8}
\def\refW     {9}
\def\refZ     {10}
\def\refZL    {11}    
\def\refZO    {12}   \def\refOZ{\refZO}

\def\C {\Bbb C}

\def\CP{\Bbb{CP}}
\def\aff{{\operatorname{aff}}}

\topmatter

\title   Remark on Tono's theorem about cuspidal curves
\endtitle

\author  S.~Yu.~Orevkov
\endauthor

\address  IMT, Universit\'e Toulouse-3, France \endaddress
\address Steklov Math. Institute, Moscow, Russia \endaddress

\endtopmatter

\document
\centerline{\hskip 45mm\it to Vladimir Lin in occasion of his 80-th birthday}

\bigskip
Lin and Zaidenberg [\refLNM, \refLNMii] (see also [\refAS; \S5], [\refLin], and Remark 2 below)
asked the following questions.

\smallskip
(Q1) {\sl Does there exist a connection between the topology of an irreducible plane affine algebraic curve
and the number of its irreducible singularities?}
(Q2) {\sl Is it true, for example, that the number of
irreducible singularities of such a curve $A$ does not exceed $1+2b_1(A)$ where $b_1(A)$ is the first
Betti number of $A$? }

\smallskip
Conjecturally, the answer to the both questions is positive.
The first and the most fascinating case of this conjecture was proven by Lin and Zaidenberg themselves
[\refZL]: if $b_1(A)=0$, then an automorphism of $\C^2$
transforms $A$ into $x^p=y^q$, in particular, $A$ has at most one singular point.
Borodzik and \D Zo\l{}\c adek [\refBZ] proved that the answer to Question (Q2)
is positive in one more particular case. Namely, if $A$ is homeomorphic to an annulus,
then $A$ has at most three singular points.

If we pass from $A$ to its closure in $\CP^2$, then the number of singular points may only
increase whereas the first Betti number may only decrease.
Thus a positive answer to (Q1) follows from the analogous conjecture for plane projective curves.
A particular case of the projective conjecture was proven in [\refOZ]:
if a projectively rigid curve in $\CP^2$ is homeomorphic to a sphere, then it
has at most 9 singular points.
Then Tono [\refTono] proved a much stronger result: if a curve in $\CP^2$ is
homeomorphic to a Riemann surface of genus $g$, then it has no more than
$(21g + 17)/2$ singular points (thus no more than $8$ when $g=0$).

The purpose of this note is to remark that Tono's arguments extend without changes
to the case of an arbitrary plane projective curve and they prove the projective
analog of the conjecture and hence, a positive answer to Question (Q1).

Let us give precise statements. Let $C$ be an algebraic curve in $\CP^2$.
A singular point of $C$ is called a {\it cusp} if $C$ has a single local analytic branch at it.
Let $s$ be the number of all singular points of $C$ and $c$ the number of cusps.
Let $b_i=b_i(C)$ be the $i$-th Betti number of $C$. So, $b_2$ is the number of irreducible
components. Let $g=g(C)$ be the total genus of $C$, i.e.,
the sum of the genera of the normalizations of all the irreducible components.

\proclaim{ Theorem 1 } If
$\bar\kappa(\CP^2\setminus C)=2$ {\rm(by [\refW] this is so, for example,
when one of irreducible components of $C$ has $\ge3$ singular points)}, then
$c \le {9\over2} b_1 + {3\over2}g - 6 b_2 + {29\over2}$ and
$s \le {11\over2} b_1 - {1\over2}g - 5 b_2 + {27\over2}$.
\endproclaim

\proclaim{ Corollary 1 } If $C$ is irreducible, then
$c \le {9\over2} b_1 + {3\over2} g + {17\over2} \le {21\over4}b_1 + {17\over2}$
and
$s \le {11\over2} b_1 - {1\over2} g + {17\over2} \le {11\over2} b_1 + {17\over2}$.
\endproclaim

Let $C^\aff$ be the intersection of $C$ with some fixed affine chart and let $b_i^\aff=b_i(C^\aff)$.
We denote the number of singular points, the number of cusps,
and the number of points at infinity of $C^\aff$ by $s^\aff$,
$c^\aff$, and $p$ respectively.

\proclaim{ Corollary 2 } If $\bar\kappa(\CP^2\setminus C)=2$, then
$c^\aff \le {9\over2}(b_1^\aff - b_0^\aff - p) + {3\over2}(g - b_2) + 19$ and
$s^\aff \le {11\over2}(b_1^\aff - b_0^\aff - p) + {1\over2}(b_2 -g)+ 19$.
\endproclaim

\proclaim{ Corollary 3 } If $C$ is irreducible, then
$c^\aff \le {9\over2}(b_1^\aff - p) + {3\over2}g + 13 \le  {9\over2} b_1^\aff + {3\over2}g + {17\over2} $ and
$s^\aff \le {11\over2}(b_1^\aff - p) - {1\over2}g + 14 \le {11\over2} b_1^\aff - {1\over2}g + {17\over2} $.
\endproclaim

\demo{ Proof of Theorem 1 }
Let $\sigma:V\to\CP^2$ be the birational morphism (a composition of blowups) such that
$D=\sum D_i=\sigma^{-1}(C)$ is a curve with simple normal crossings all whose irreducible components
$D_i$ are smooth. Let $K$ be the canonical class of $V$ and let $p_a(D) = D(K + D)/2 + 1$.
Then [\refTono; Corollary 4.4] combined with the lower bound $2c$ for the number of maximal twigs
(proved by the same arguments as in [\refTono; p.~220]) yields
$$
      2c \le 12 e(V\setminus D) + 5 - 3 p_a(D)             \eqno(1)
$$
where $e$ stands for the Euler characteristic.
Let $\Gamma$ be the dual graph of $D$, i.e., the vertices of $\Gamma$ correspond to the irreducible
components of $D$ and the edges correspond to the crossing points. Let $b_1^\Gamma=b_1(\Gamma)$.
Then we have $p_a(D) = g + b_1^\Gamma$ and $b_1 =  b_1(D) = 2g + b_1^\Gamma$
whence $p_a(D) = b_1-g$. We have also
$e(V\setminus D) = e(\CP^2\setminus C) = 2 + b_1 - b_2$. Thus (1) yields the required
bound for $c$.

Let $d=\sum_{i=1}^s (r_i-1)$ where
$r_i$ is the number of local analytic branches of $C$ at the $i$-th singular point.
We have $d \ge s-c$.
We define the {\it incidence graph\/} of a curve as the bipartite graph whose vertices correspond
to singular points and to irreducible components of the curve, and the edges correspond to its local branches
at singular points: the edge corresponding to a local branch at $p$ connects
the vertex corresponding to $p$ with the vertex corresponding
to the irreducible component containing the local branch.
Let $\Gamma_C$ and $\Gamma_D$ be the incidence graphs of $C$ and $D$.
It is clear that $\Gamma_D$ is homeomorphic to $\Gamma$. Since the
topology of incidence graphs does not change under blowups,
it follows that $\Gamma_C\cong\Gamma_D\cong\Gamma$.
Further, $\Gamma_C$ has $b_2+s$ vertices and $d+s$ edges, thus
$e(\Gamma_C) = b_2 - d$. Since $e(\Gamma) = 1 - b_1^\Gamma = 1 - b_1+2g$,
we obtain 
$$
   s \le c + d = c + b_2 - e(\Gamma) = c + b_2 + b_1 - 2g - 1
$$
and the result follows from the bound for $c$.
\enddemo

\demo{ Proof of Corollaries }
{\bf 1.} Since the upper bounds for $c$ and $s$ are $\ge3$, it is enough to
consider the case when $s\ge3$ and hence $\bar\kappa(\CP^2\setminus C)\ge2$ by [\refW].

{\bf 2, 3.} Use $c^\aff\le c$, $s^\aff\le s$, and
$$
    b_1 = b_2+1-e(C) = b_2+1-e(C^\aff)-p = b_2+1-b_0^\aff +b_1^\aff-p
$$
in Theorem 1 and Corollary 1.
\enddemo

\smallskip\noindent
{\bf Remark 1.} In a recent preprint [\refPalka; Theorem 1.4], Palka improved Tono's estimate
$c\le 8$ in the case of rational cuspidal curves up to $c\le 6$. Maybe,
his arguments could give a better upper bound for $c$ and $c^\aff$ in the general case.

\smallskip\noindent
{\bf Remark 2.} In the case of an irreducible affine curve, Corollary 3 gives an estimate
$c^\aff\le\alpha b_1^\aff+\beta$ with $(\alpha,\beta)=(5.25,8.5)$. This estimate does not seem
to be optimal. The example $y^2=p(x)^3$ shows that one cannot do better than $(\alpha,\beta)=(1,1)$ (see also [\refAS]).
Note that the projective dual of a smooth or nodal cubic with a suitable choice of the infinite line
are the only known exceptions for the estimate with $(\alpha,\beta)=(1,1)$; for them we have
$b_1^\aff=1,5,6$, or $7$.

\smallskip\noindent
{\bf Remark 3.} Let $D=D_1+\dots+D_n$ be a reduced curve with simple normal crossings
on a smooth algebraic surface $V$ and let $\Gamma$ be the dual graph of $D$.
We set $\beta(D_i)=D_i(D-D_i)$ (the degree of the corresponding vertex of $\Gamma$). If $\beta(D_i)=1$,
we say that $D_i$ is a {\it tip} of $D$. 
We assume that $D$ does not contain any rational $(-1)$-curve
$D_i$ with $\beta(D_i)\le 2$.
Zaidenberg [\refZ; p.~16] conjectured that {\sl only finite number of pairwise non-homeomorphic graphs $\Gamma$
can be obtained in this way under the condition that $\bar\kappa(V\setminus D)=2$ and
$b_i(V\setminus D)=0$ for $i>0$.}

Similarly to Theorem 1, this conjecture follows immediately from [\refTono; Corollary 4.4]. Indeed,
the condition $b_i(V\setminus D)=0$, $i>0$, implies that
$\Gamma$ is a tree. Hence it is enough to bound $t_D$ (the number of tips of $D$).
A bound $t_D\le 17$ follows from [\refTono; Corollary 4.4] provided that the pair $(V,D)$ is
{\it almost minimal} which in our setting is equivalent to the absence of rational $(-1)$-curves $E$ such that
$E\not\subset D$ and $ED\le 1$ (see [\refTono; Lemma 3.4]). If such $E$ exists, let $\pi:V\to V'$ be
its blowing down and let $D'=\pi_*(D)$. Then $\bar\kappa(V'\setminus D')=2$ by [\refTono; Lemma 3.2]
and $e(V'\setminus D')\le e(V\setminus D)-1=0$. This contradicts [\refM; Theorem 6.7.1].

There are 3901520 trees with at most 17 tips (the one vertex tree included).

\Refs
\def\r{\ref}

\r\no\refAS
\by     S.S.~Abhyankar and A.~Sathaye
\paper  Uniqueness of plane embeddings of special curves
\jour   Proc. Amer. Math. Soc. \vol 124 \yr 1996 \pages 1061--1069
\endref

\r\no\refBZ
\by     M.~Borodzik and H.~\D Zo\l{}\c adek
\paper  Number of singular points of an annulus in $\Bbb C^2$
\jour   Ann. Inst. Fourier \vol 61 \yr 2011 \pages 1539--1555
\endref

\r\no\refLin
\by     V.~Lin
\paper  Some problems that I would like to see solved
\jour   Abstract of a talk. Technion, 2015,
http://www2.math.technion.ac.il/$\widetilde{\;}$pincho/Lin/Abstracts.pdf
\endref

\r\no\refLNM
\by       V.Ya.~Lin and M.G.~Zaidenberg
\paper    On the number of singular points of a plane affine algebraic curve
\inbook   in: Linear and complex analysis problem book
\bookinfo Lect. Notes Math. \vol 1043 \publ Springer \yr 1984 \pages 662--663
\endref

\r\no\refLNMii
\by       V.Ya.~Lin and M.G.~Zaidenberg
\paper    On the number of singular points of a plane affine algebraic curve
\inbook   in: Linear and complex analysis problem book
\bookinfo Lect. Notes Math. \vol 1574 \publ Springer \yr 1994 \page 479
\endref

\r\no\refM
\by       M.~Miyanishi
\book     Open algebraic surfaces
\publ     A. M. S.
\publaddr Providence RI \yr 2001
\endref

\r\no\refPalka
\by     K.~Palka
\paper  Cuspidal curves, minimal models and Zaidenberg's finiteness conjecture
\jour \hbox to 10mm{}\ arXiv:1405.5346
\endref

\r\no\refTono
\by     K.~Tono
\paper  On the number of the cusps of cuspidal plane curves
\jour   Math. Nachr. \vol 278 \yr 2005 \pages 216--221
\endref

\r\no\refW
\by     I.~Wakabayashi
\paper  On the logarithmic Kodaira dimension of the complement of a curve in $P^2$
\jour   Proc. Japan Acad. Ser A \vol 54 \yr 1978 \pages 157--162
\endref

\r\no\refZ
\by    M.~Zaidenberg
\paper in: Open problems on open algebraic varieties \jour arXiv:alg-geom/9506006
\endref

\r\no\refZL
\by    M.G.~Zaidenberg and V.Ya.~Lin
\paper An irreducible, simply connected algebraic curve in $\Bbb C^2$
       is equivalent to a quasihomogeneous curve
\jour  Dokl. Akad. Nauk SSSR \vol 271:5 \yr 1983 \pages 1048--1052
\lang  Russian \transl English Transl.
\jour  Soviet Math., Doklady \vol 28 \yr 1983 \pages 200--204
\endref

\r\no\refZO
\by    M.G.~Zaidenberg and S.Yu.~Orevkov
\paper On rigid rational cuspidal plane curves
\jour  Uspekhi Mat. Nauk \vol 51:1 \yr 1996 \pages 149--150
\lang  Russian \transl English Transl.
\jour  Russian Math. Surveys \vol 51 \yr 1996 \pages 179--180
\endref

\endRefs
\enddocument